\documentclass[11pt]{amsart}
\newtheorem{theorem}{Theorem}[section]

\newtheorem{lemma}[theorem]{Lemma}

\newtheorem{example}[theorem]{Example}
\theoremstyle{definition}

\theoremstyle{remark}
\newtheorem{remark}{Remark}
\numberwithin{equation}{section}
\begin{document}

\title[Non-energy semi-stable radial solutions]
{Non-energy semi-stable radial solutions}
\author{Salvador Villegas}
\thanks{The author has been supported by the MEC Spanish grant MTM2012-37960.}
\address{Departamento de An\'{a}lisis
Matem\'{a}tico, Universidad de Granada, 18071 Granada, Spain.}
\email{svillega@ugr.es}

\begin{abstract} This paper is devoted to the study of semi-stable radial solutions $u\notin H^1(B_1)$ of $-\Delta u=f(u) \mbox{ in } \overline{B_1}\setminus \{ 0\}=\{ x\in \mathbb{R}^N : 0<\vert x\vert\leq 1\}$, where $f\in C^1(\mathbb{R})$ and $N\geq 2$. We establish sharp pointwise estimates for such solutions. In addition, we prove that in dimension $N=2$, any semi-stable radial weak solution of $-\Delta u=f(u)$, posed in $B_1$ with Dirichlet data $u|_{\partial B_1}=0$, is regular.
\end{abstract}

\maketitle
\section{Introduction and main results}

This paper deals with the semi-stability of radial solutions of

\begin{equation}\label{mainequation}
-\Delta u=f(u) \ \ \mbox{ in } \overline{B_1}\setminus \{ 0\},
\end{equation}

\noindent where $B_1$ is the open unit ball of $\mathbb{R}^N$, $N\geq 2$ and $f\in C^1(\mathbb{R})$. We consider classical solutions $u\in C^2(\overline{B_1}\setminus \{ 0\})$. This is not a restriction. In fact, if we consider a radial solution $u$ of this equation in a very weak sense, we obtain that $u$ is a $C^3$ function.

A solution $u$ of (\ref{mainequation}) is called semi-stable if

$$Q_u(v):=\int_{B_1} \left( \vert \nabla v\vert^2-f'(u)v^2\right) \, dx\geq 0$$

\noindent for every $v\in C^1 (B_1)$ with compact support in $B_1\setminus\{ 0\}$. Formally, the above expression is the second variation of the energy functional associated to (\ref{mainequation}) in a domain $\Omega\subset\subset B_1\setminus \{ 0\}$: $E_\Omega (u)=\int_\Omega \left( \vert \nabla u\vert^2 /2-F(u)\right) \, dx$, where $F'=f$. Thus, if $u\in C^1(\overline{B_1}\setminus \{0 \})$ is a local minimizer of $E_\Omega$ for every smooth domain $\Omega\subset\subset B_1\setminus \{ 0\}$ (i.e. a minimizer under every small enough $C^1(\overline{\Omega})$ perturbation vanishing on $\partial \Omega$), then $u$ is a semi-stable solution of (\ref{mainequation}).

We will be also interested in the semi-stability of radial weak solutions of the problem
\begin{equation}\label{theball}
\left\{
\begin{array}{ll}
-\Delta u=f(u)&\mbox{in } B_1 \, ,\\
u=0&\mbox{on } \partial B_1 \, ,\\
\end{array}
\right.
\end{equation}

\noindent where  $N\geq 2$ and $f\in C(\mathbb{R})$.

As in \cite{Bre}, we say that $u$ is a weak solution of (\ref{theball}) if $u\in L^1(B_1)$, $f(u)\delta\in L^1(B_1)$ and

\begin{equation}\label{weakdefinition}
-\int_{B_1}u\Delta \zeta dx=\int_{B_1}f(u)\zeta dx
\end{equation}

\noindent for all $\zeta\in C^2(\overline{B_1})$ with $\zeta=0$ on $\partial B_1$. Here $\delta (x)=\mbox{dist}(x,\partial B_1)$ denotes the distance to the boundary of $B_1$.

If $f\in C^1(\mathbb{R})$, we say that a radial weak solution $u$ of (\ref{theball}) is semi-stable if $u|_{\overline{B_1}\setminus\{ 0\} }$ is semi-stable. This definition has sense, since any radial weak solution of (\ref{theball}) is a $C^2(\overline{B_1}\setminus\{ 0\} )$ function (see Lemma \ref{characterization} below).

The original motivation of this work is the following. Consider the semilinear elliptic problem

$$
\left\{
\begin{array}{ll}
-\Delta u=\lambda g(u)\ \ \ \ \ \ \  & \mbox{ in } \Omega \, ,\\
u\geq 0 & \mbox{ in } \Omega \, ,\\
u=0  & \mbox{ on } \partial\Omega \, ,\\
\end{array}
\right. \eqno{(P_\lambda)}
$$
\

\noindent where $\Omega\subset\mathbb{R}^N$ is a smooth bounded domain, $N\geq 1$, $\lambda\geq 0$ is a real parameter and the nonlinearity $g:[0,\infty)\rightarrow \mathbb{R}$ satisfies

\begin{equation}\label{convexa}
g \mbox{ is } C^1, \mbox{ nondecreasing and convex, }g(0)>0,\mbox{
and }\lim_{u\to +\infty}\frac{g(u)}{u}=+\infty.
\end{equation}

It is well known that there exists a finite positive extremal parameter $\lambda^\ast$ such that ($P_\lambda$) has a minimal classical solution $u_\lambda\in C^2(\overline{\Omega})$ if $0\leq \lambda <\lambda^\ast$, while no solution exists, even in the weak sense (similar definition as the case $\Omega=B_1$), for $\lambda>\lambda^\ast$. The set $\{u_\lambda:\, 0\leq \lambda < \lambda^\ast\}$ forms a branch of classical solutions increasing in $\lambda$. Its increasing pointwise limit $u^\ast(x):=\lim_{\lambda\uparrow\lambda^\ast}u_\lambda(x)$ is a weak solution of ($P_\lambda$) for $\lambda=\lambda^\ast$, which is called the extremal solution of ($P_\lambda$) (see \cite{Bre,BV,Dup}).

The regularity and properties of extremal solutions depend strongly on the dimension $N$, domain $\Omega$ and nonlinearity $g$. When $g(u)=e^u$, it is known that $u^\ast\in L^\infty (\Omega)$ if $N<10$ (for every $\Omega$) (see \cite{CrR,MP}), while $u^\ast (x)=-2\log \vert x\vert$ and $\lambda^\ast=2(N-2)$ if $N\geq 10$ and $\Omega=B_1$ (see \cite{JL}). There is an analogous result for $g(u)=(1+u)^p$ with $p>1$ (see \cite{BV}). Brezis and V\'azquez \cite{BV} raised the question of determining the boundedness of $u^\ast$, depending on the dimension $N$, for general nonlinearities $g$ satisfying (\ref{convexa}). The first general results were due to Nedev \cite{Ne}, who proved that $u^\ast \in L^\infty (\Omega)$ if $N\leq 3$, and $u^\ast \in L^p (\Omega)$ for every $p<N/(N-4)$, if $N\geq 4$. In a recent paper the author \cite{yo} has proved that $u^\ast \in L^\infty (\Omega)$ if $N=4$, and $u^\ast \in L^{N/(N-4)}(\Omega)$, if $N\geq 5$. Cabr\'e \cite{cabre4}, proved that $u^\ast \in L^\infty (\Omega)$ if $N\leq 4$ and $\Omega$ is convex (no convexity on $f$ is imposed). If $N\geq 5$ and $\Omega$ is convex Cabr\'e and Sanch\'on \cite{casa} have obtained that $u^\ast \in L^\frac{2N}{N-4} (\Omega)$ (again, no convexity on $f$ is imposed). On the other hand, Cabr\'e and Capella \cite{cabrecapella} have proved that $u^\ast \in L^\infty (\Omega)$ if $N\leq 9$ and $\Omega=B_1$. Recently, Cabr\'e and Ros-Oton \cite{cros} have obtained that $u^\ast \in L^\infty (\Omega)$ if $N\leq 7$ and $\Omega$ is a convex domain of double revolution (see \cite{cros} for the definition).

Another interesting question is whether the extremal solution lies in the energy class.  Nedev \cite{Ne,Ne2} proved that $u^\ast \in H_0^1(\Omega)$ if $N\leq 5$ (for every $\Omega$) or $\Omega$ is convex (for every $N\geq 1$). The author \cite{yo} has obtained that $u^\ast \in H_0^1(\Omega)$ if $N=6$ (for every $\Omega$). Brezis and V\'azquez \cite{BV} proved that a sufficient condition to have $u^\ast \in H_0^1(\Omega)$ is that $\liminf_{u\to \infty} u\, g'(u)/g(u)>1$ (for every $\Omega$ and $N\geq 1$).

Note that the minimality of $u_\lambda$ ($0<\lambda<\lambda^\ast$) implies its semi-stability, i.e. $\int_\Omega \left( \vert \nabla v\vert^2-\lambda g'(u_\lambda)v^2\right) \, dx\geq 0$, for every $v\in C^1(\Omega)$ with compact support. Clearly, we can pass to the limit and obtain that $u^\ast$ is also a semi-stable weak solution for $\lambda=\lambda^\ast$. Conversely, in \cite{BV} it is proved that if $g$ satisfies (\ref{convexa}) and $u\in H_0^1(\Omega)$ is an unbounded semi-stable weak solution of ($P_\lambda$) for some $\lambda>0$, then $u=u^\ast$ and $\lambda=\lambda^\ast$. (For instance, applying this result it follows easily that $u^\ast (x)=-2\log \vert x\vert$ and $\lambda^\ast=2(N-2)$ if $g(u)=e^u$, $\Omega=B_1$ and $N\geq 10$). The hypothesis $u\in H_0^1(\Omega)$ is essential since in \cite{BV} it is observed that if $\Omega=B_1$, $N\geq 3$ and  $\displaystyle{\frac{N}{N-2}<p\leq \frac{N+2\sqrt{N-1}}{N+2\sqrt{N-1}-4}}$, then $u(x)=\vert x\vert^{-2/(p-1)}-1$ is an unbounded semi-stable weak solution of ($P_\lambda$) for $g(u)=(1+u)^p$ and $\lambda=2(Np-2p-N)/(p-1)^2$, which is a non-energy function, i.e. $u\notin H_0^1(B_1)$. Since $B_1$ is a convex domain, $u^\ast\in H_0^1(B_1)$ and then $u\not\equiv u^\ast$. As pointed out in \cite{BV}, this type of "strange" solutions are apparently isolated objects that cannot be obtained as limit of classical solutions, which leaves them in a kind of "limbo" with respect to the classical theory.

In this paper we study this class of non-energy semi-stable radial solutions and it is established sharp pointwise estimates for such solutions. In addition we prove that, contrary to the case $N\geq 3$, there is no solutions of this type in dimension $N=2$.

\begin{theorem}\label{nonenergy}

Let $N\geq 2$, $f\in C^1(\mathbb{R})$ and $u\notin H^1(B_1)$ be a semi-stable radial solution of (\ref{mainequation}). Then there exist $M>0$ and $0<r_0<1$ such that

$$\vert u(r)\vert\geq\left\{
\begin{array}{ll}
M \vert\log r\vert & \forall r\in (0,r_0) \mbox{ if } N=2, \\ \\
M r^{-N/2-\sqrt{N-1}+2} &  \forall r\in (0,r_0) \mbox{ if } N\geq 3. \\
\end{array}
\right.
$$

\end{theorem}

\

\begin{theorem}\label{positivef}

Let $N\geq 2$, $0\leq f\in C^1(\mathbb{R})$ and $u\notin H^1(B_1)$ be a semi-stable radially decreasing near the origin solution of (\ref{mainequation}). We have that:

\begin{enumerate}

\item[i)] If $N=2$, then $\lim_{r\to 0} ru'(r)=-\alpha$, for some $\alpha\in (0,+\infty)$. In particular $\lim_{r\to 0} u(r)/\vert \log r\vert =\alpha$.

\

\item[ii)] If $N\geq 3$, then $M_1 r^{-N/2-\sqrt{N-1}+1}\leq \vert u_r(r)\vert \leq M_2 r^{-N+1}$ in $\overline{B_1}$, for some constants $M_1, M_2 >0$.

\end{enumerate}

\end{theorem}

\

\begin{theorem}\label{N=2}

Let $N=2$, $f\in C^1(\mathbb{R})$ and $u$ be a semi-stable radial weak solution of (\ref{theball}). Then $u$ is regular (i.e. $u\in C^2(\overline{B_1})$).

\end{theorem}

The main results obtained in this paper are optimal. If $N=2$, clearly $u(r)=\vert \log r\vert\notin H^1(B_1)$ satisfies $-\Delta u=0$ and then it is a semi-stable radial solution of (\ref{mainequation}) for $f\equiv 0$. On the other hand, for every $N\geq 2$ and $\alpha<0$ consider the radial function $u_\alpha(r)=r^\alpha$, $0<r\leq 1$ and a function $f_\alpha\in C^\infty(\mathbb{R})$ satisfying $f_\alpha(s)=-\alpha\left(\alpha+N-2\right)s^{1-2/\alpha}$ for every $s\geq 1$. If $N\geq 3$ and $2-N\leq \alpha<0$ then we take $f_\alpha\geq 0$. The following example shows that the pointwise estimates of Theorems \ref{nonenergy} and \ref{positivef} are sharp.

\begin{example}
Let $\alpha<0$ if $N=2$ and $\alpha\leq-N/2-\sqrt{N-1}+2$ if $N\geq 3$. Consider the above defined functions $u_\alpha, f_\alpha$. Then $u_\alpha\notin H^1(B_1)$ is a semi-stable radial solution of (\ref{mainequation}) for $f=f_\alpha$.
\end{example}

\noindent\textbf{Proof.}
It is immediate that $u_\alpha\notin H^1(B_1)$ is a radial solution of (\ref{mainequation}) for $f=f_\alpha$. An easy computation shows that $f'_\alpha(u_\alpha(r))=-(\alpha -2)(\alpha+N-2)/r^2$, for every $0<r\leq 1$. Taking into account that $\alpha<0$ if $N=2$, and $\alpha\leq -N/2-\sqrt{N-1}+2$ if $N=3$, we check at once that $-(\alpha -2)(\alpha+N-2)\leq (N-2)^2/4$, which is the best constant in Hardy's inequality: $\int_{B_1} ((N-2)^2/(4r^2)) v^2\leq \int_{B_1} \vert \nabla v \vert^2$, for every $v\in C^1 (B_1)$ with compact support in $B_1\setminus\{ 0\}$. This gives the semi-stability of $u_\alpha$ for this range of values of $\alpha$. \qed

\section{Sharp pointwise estimates}

Lemmas \ref{maintool} and \ref{monotonicity} below are almost identical to Lemmas 2.1 and 2.2 of \cite{yomismo}. We prefer to state them here and give the same proof as in \cite{yomismo} for the convenience of the reader. In fact, Lemma \ref{maintool} follows easily from the ideas of the proof of \cite[Lem. 2.1]{cabrecapella}, which was inspired by the proof of Simons theorem on the nonexistence of singular minimal cones in $\mathbb{R}^N$ for $N\leq 7$ (see \cite[Th. 10.10]{simons} and \cite[Rem. 2.2]{cabrecapella} for more details).

\begin{lemma}\label{maintool}

Let $N\geq 2$, $f\in C^1(\mathbb{R})$ and $u$ be a semi-stable radial solution of (\ref{mainequation}).  Let $0<r_1<r_2<1$ and $\eta\in C^{0,1}([r_1,r_2])$ such that $\eta u_r$ vanishes at $r=r_1$ and $r=r_2$. Then

$$\int_{r_1}^{r_2} r^{N-1}u_r^2 \left( \eta'^2-\frac{N-1}{r^2}\eta^2\right) dr\geq 0.$$

\end{lemma}

\noindent\textbf{Proof.}
First of all, note that we can extend the second variation of energy $Q_u$ to the set of functions $v\in C^{0,1}(B_1)$ with compact support in $B_1\setminus\{ 0\}$, obtaining $Q_u(v)\geq 0$ for such functions $v$. Hence, we can take the radial function $v=\eta u_r \chi_{B_{r_2}\setminus\overline{B_{r_1}}}$.

On the other hand, differentiating (\ref{mainequation}) with respect to r, we have

$$-\Delta u_r+\frac{N-1}{r^2}u_r=f'(u)u_r, \ \ \mbox{ for all }r\in (0,1).$$

Following the ideas of the proof of \cite[Lem. 2.1]{cabrecapella}, we can multiply this equality by $\eta^2 u_r$ and integrate by parts in the annulus of radii $r_1$ and $r_2$ to obtain

$$0=\int_{B_{r_2}\setminus\overline{B_{r_1}}}\left(\nabla u_r\nabla\left(\eta^2 u_r\right)+\frac{N-1}{r^2}u_r \eta^2 u_r-f'(u)u_r\eta^2 u_r\right)dx$$
$$=\int_{B_{r_2}\setminus\overline{B_{r_1}}}\left(\vert \nabla\left(\eta u_r\right)\vert^2-f'(u)\left( \eta u_r\right)^2\right)dx-\int_{B_{r_2}\setminus\overline{B_{r_1}}} u_r^2\left(\vert\nabla\eta\vert^2-\frac{N-1}{r^2}\eta^2\right)dx$$
$$=Q_u(\eta u_r \chi_{B_{r_2}\setminus\overline{B_{r_1}}})-\omega_N\int_{r_1}^{r_2} r^{N-1}u_r^2 \left( \eta'^2-\frac{N-1}{r^2}\eta^2\right) dr.$$

Using the semi-stability of $u$ the lemma follows. \qed

\begin{lemma}\label{monotonicity}

Let $N\geq 2$, $f\in C^1(\mathbb{R})$ and $u$ be a nonconstant semi-stable radial solution of (\ref{mainequation}). Then $u_r$ vanishes at most in one value in $(0,1)$.

\end{lemma}

\noindent\textbf{Proof.}
Suppose by contradiction that there exist $0<r_1<r_2<1$ such that $u_r(r_1)=u_r(r_2)=0$. Taking $\eta\equiv 1$ in the previous lemma, we obtain

$$\int_{r_1}^{r_2} r^{N-1}u_r^2 \left( -\, \frac{N-1}{r^2}\right) dr\geq 0.$$

Hence we conclude that $u_r\equiv 0$ in $[r_1,r_2]$, which clearly forces $u$ is constant in $\overline{B_1}\setminus\{ 0\}$, a contradiction. \qed

\begin{lemma}\label{theidea}

Let $N\geq 2$, $u\in C^2(\overline{B_1}\setminus \{ 0\})$ a radial function satisfying $u\notin H^1(B_1)$. Then there exist $0<a<1/2$ and a function $\eta_0\in C^{0,1}([a,1/2])$ such that $\eta_0(a)=1$, $\eta_0(1/2)=0$ and

$$\int_{a}^{1/2} r^{N-1}u_r^2 \left( \eta_0'^2-\frac{N-1}{r^2}\eta_0^2\right) dr<0.$$ \qed

\end{lemma}

\noindent\textbf{Proof.}
For arbitrary $a\in (0,1/4)$ define the function

$$\eta_0 (r)=\left\{
\begin{array}{ll}
1 & \mbox{ if } a\leq r <1/4, \\ \\
2-4r & \mbox{ if } 1/4\leq r\leq 1/2. \\
\end{array}
\right.
$$

Clearly $\eta_0$ is a $C^{0,1}([a,1/2])$ function satisfying $\eta_0(a)=1$ and $\eta_0(1/2)=0$. On the other hand

$$\int_{a}^{1/2} r^{N-1}u_r^2 \left( \eta_0'^2-\frac{N-1}{r^2}\eta_0^2\right) dr=-(N-1)\int_{a}^{1/4}r^{N-3}u_r^2\, dr$$ $$+ \int_{1/4}^{1/2} r^{N-1}u_r^2 \left( 16-\frac{N-1}{r^2}(2-4r)^2\right) dr.$$

Note that $u\in C^2(\overline{B_1}\setminus \{ 0\})$ and $u\notin H^1(B_1)$ imply $r^{N-1}u_r^2\notin L^1(0,1/4)$ and therefore $\int_0^{1/4}r^{N-3}u_r^2\, dr=+\infty$. From the above it follows that

$$\lim_{a\to 0}\int_{a}^{1/2} r^{N-1}u_r^2 \left( \eta_0'^2-\frac{N-1}{r^2}\eta_0^2\right) dr=-\infty.$$

Taking $a\in (0,1/4)$ sufficiently small the lemma follows. \qed

\begin{lemma}\label{largelemma1}
Let $N\geq 2$, $f\in C^1(\mathbb{R})$ and $u\notin H^1(B_1)$ be a semi-stable radial solution of (\ref{mainequation}). Then there exist $K>0$ and $0<r_0<1$ such that

$$\int_{r/2}^r\frac{ds}{u_r(s)^2}\leq K r^{N+2\sqrt{N-1}-1}\, \ \ \ \forall r\in (0,r_0).$$

\end{lemma}

\

\noindent\textbf{Proof.}
Consider $a$ and $\eta_0$ of Lemma \ref{theidea}. From Lemma \ref{monotonicity} we can choose $0<r_0<a$ such that $u_r$ does not vanish in $(0,r_0]$. We now fix $r\in (0,r_0)$ and consider the function

$$\eta (t)=\left\{
\begin{array}{ll}
\displaystyle{\frac{r^{\sqrt{N-1}}}{\int_{r/2}^r
\frac{ds}{u_r(s)^2}}} \int_{r/2}^t \frac{ds}{u_r(s)^2} &
\mbox{ if } r/2\leq t\leq r, \\ \\
t^{\sqrt{N-1}} & \mbox{ if } r<t\leq a, \\ \\
a^{\sqrt{N-1}}\eta_0(t) & \mbox{ if } a <t\leq 1/2. \\
\end{array}
\right.
$$

\

Applying Lemma \ref{maintool} (with $r_1=r/2$ and $r_2=1/2$) we obtain

$$0\leq \int_{r/2}^{1/2} t^{N-1}u_r(t)^2 \left( \eta'(t)^2-\frac{N-1}{t^2}\eta(t)^2\right) dt$$
$$=\int_{r/2}^r t^{N-1}u_r(t)^2 \left( \eta'(t)^2-\frac{N-1}{t^2}\eta(t)^2\right) dt$$
$$+ a^{2\sqrt{N-1}}\int_{a}^{1/2} t^{N-1}u_r(t)^2 \left( \eta_0'(t)^2-\frac{N-1}{t^2}\eta_0(t)^2\right) dt$$
$$\leq r^{N-1}\int_{r/2}^r u_r(t)^2 \eta'(t)^2 dt+a^{2\sqrt{N-1}}\int_{a}^{1/2} t^{N-1}u_r(t)^2 \left( \eta_0'(t)^2-\frac{N-1}{t^2}\eta_0(t)^2\right) dt$$
$$=r^{N-1}\frac{r^{2\sqrt{N-1}}}{\int_{r/2}^r\frac{ds}{u_r(s)^2}}+a^{2\sqrt{N-1}}\int_{a}^{1/2} t^{N-1}u_r(t)^2 \left( \eta_0'(t)^2-\frac{N-1}{t^2}\eta_0(t)^2\right) dt.$$

This gives

$$-a^{2\sqrt{N-1}}\int_{a}^{1/2} t^{N-1}u_r(t)^2 \left( \eta_0'(t)^2-\frac{N-1}{t^2}\eta_0(t)^2\right) dt\leq \displaystyle{\frac{r^{N+2\sqrt{N-1}-1}}{\int_{r/2}^r\frac{ds}{u_r(s)^2}}},$$

\noindent which is the desired conclusion for

$$K=\left( -a^{2\sqrt{N-1}}\int_{a}^{1/2} t^{N-1}u_r(t)^2 \left( \eta_0'(t)^2-\frac{N-1}{t^2}\eta_0(t)^2\right) dt\right) ^{-1},$$

\noindent which is a positive number, from Lemma \ref{theidea}. \qed

\begin{lemma}\label{largelemma2}

Let $N\geq 2$, $f\in C^1(\mathbb{R})$ and $u\notin H^1(B_1)$ be a semi-stable radial solution of (\ref{mainequation}). Then there exist $M'>0$ and $0<r_0<1$ such that

$$\vert u(r)-u(r/2)\vert \geq M' r^{-N/2-\sqrt{N-1}+2}\,  \ \ \ \forall r\in (0,r_0).$$
\end{lemma}

\

\noindent\textbf{Proof.}
Take the same constant $0<r_0<1$ of Lemma \ref{largelemma1}. Fix $r\in (0,r_0)$ and consider the functions:

$$\begin{array}{ll}
 \alpha (s)=\vert u_r(s)\vert^{-\frac{2}{3}}, & s\in (r/2,r). \\
 \\
  \beta (s)=\vert u_r(s)\vert^{\frac{2}{3}}, & s\in (r/2,r). \\

\end{array}
$$

By Lemma \ref{largelemma1} we have

$$\Vert \alpha \Vert _{L^3(r/2,r)}\leq K^\frac{1}{3} \,
r^\frac{N+2\sqrt{N-1}-1}{3}$$

\noindent for a constant $K>0$ not depending on $r\in (0,r_0)$. On the other hand, since $u_r$ does not vanish in $(0,a]$, it follows

$$\Vert \beta \Vert _{L^{3/2}(r/2,r)}=\vert u(r)-u(r/2)\vert^\frac{2}{3}.$$

Applying H{\"o}lder inequality to functions $\alpha$ and $\beta$ we deduce

$$r/2=\int_{r/2}^r \alpha(s)\beta(s) ds\leq \Vert \alpha \Vert _{L^3(r/2,r)} \Vert \beta \Vert _{L^{3/2}(r/2,r)}\leq K^\frac{1}{3}r^\frac{N+2\sqrt{N-1}-1}{3}\, \vert u(r)-u(r/2) \vert^\frac{2}{3},$$

\noindent which is the desired conclusion for $M'=2^{-3/2}K^{-1/2}$. \qed

\

\noindent\textbf{Proof of Theorem \ref{nonenergy}.}
Consider the numbers $M'>0$ and $0<r_0<1$ of Lemma \ref{largelemma2}. It is easily seen that for every $r\in (0,r_0)$ there exist an integer $m\geq 0$ and $r_0/2\leq z<r_0$ such that $r=z/2^m$. From the monotonicity of $u$ in $(0,r_0)$ it follows that

\begin{equation}\label{ww}
\vert u(r)\vert \geq \vert u(z)-u(r)\vert -\vert u(z)\vert
=\sum_{k=0}^{m-1}\left\vert u\left(\frac{z}{2^k}\right)-u\left(\frac{z}{2^{k+1}}\right)\right\vert \, -\vert
u(z)\vert
\end{equation}

$\bullet$ \textbf{Case $N=2$.} We have that $-N/2-\sqrt{N-1}+2=0$. Hence, applying Lemma \ref{largelemma2} and (\ref{ww}) we obtain

$$\vert u(r)\vert \geq M'm -\vert u(z)\vert=\frac{M'(\log z -\log r)}{\log 2}-\vert u(z)\vert ,$$

\noindent where $M'>0$ does not depend on $r\in (0,r_0)$. Since $z\in [r_0/2,r_0)$ and $u$ is continuous the above inequality is of the type

$$\vert u(r)\vert \geq M_1 \vert \log r\vert-M_2 \, \ \ \ \ \
\ \forall r\in (0,r_0),$$

\noindent for certain $M_1, M_2>0$. Taking a smaller $0<r_0<1$ if necessary, the theorem is proved in this case.

\

$\bullet$ \textbf{Case $N\geq 3$.} We have that $-N/2-\sqrt{N-1}+2<0$. Thus, applying again Lemma \ref{largelemma2} and (\ref{ww}) we deduce

$$\vert u(r)\vert \geq \sum_{k=0}^{m-1}M'\left(\frac{z}{2^k}\right)^{-N/2-\sqrt{N-1}+2} \, -\vert u(z)\vert$$ $$=M'\left(\frac{
r^{-N/2-\sqrt{N-1}+2}-z^{-N/2-\sqrt{N-1}+2}}{2^{N/2+\sqrt{N-1}-2}-1}\right)-\vert u(z)\vert ,$$

\noindent which is an inequality of the type $\vert u(r)\vert \geq M_1 r^{-N/2-\sqrt{N-1}+2}-M_2, \forall r\in (0,r_0)$, for certain $M_1, M_2>0$.  The proof is complete as the previous case. \qed

\

\noindent\textbf{Proof of Theorem \ref{positivef}.}

If $N=2$, then $(-ru_r(r))'=rf(u(r))\geq 0$ for every $r\in (0,1]$. Since $-r u_r(r)$ is nonnegative for small $r$, it is deduced that $\lim_{r\to 0}(- ru_r(r))=\alpha$, for some $\alpha \in [0,\infty)$. This implies $\lim_{r\to 0} u(r)/\vert \log r\vert=\alpha$. Applying Theorem \ref{nonenergy} we deduced $\alpha>0$, which is our claim for $N=2$.

If $N\geq 3$, then $(-r^{N-1}u_r(r))'=r^{N-1}f(u(r))\geq 0$ for every $r\in (0,1]$. Since $-r^{N-1}u_r(r)$ is nonnegative for small $r$, it is deduced that $-r^{N-1}u_r(r)$ is a nonnegative nondecreasing function and then $r^{N-1}\vert u_r(r)\vert=-r^{N-1}u_r(r)\leq -u_r(1)$, following the second inequality of (ii) for $M_2=-u_r(1)$. (Note that we have used neither the semi-stability of $u$ nor $u\notin H^1(B_1)$). To prove the first inequality of (ii), let us observe that since $-r^{N-1}u_r$ is a nonnegative nondecreasing function then $r^{2N-2}u_r^2$ is nondecreasing. Then applying Lemma \ref{largelemma1} we have that there exist $K>0$ and $0<r_0<1$ such that

$$K r^{N+2\sqrt{N-1}-1}\geq \int_{r/2}^r\frac{ds}{u_r(s)^2}=\int_{r/2}^r \frac{s^{2N-2}}{s^{2N-2}u_r(s)^2}ds$$ $$\geq\frac{1}{r^{2N-2}u_r(r)^2}\int_{r/2}^r s^{2N-2}ds=\frac{(1-2^{1-2N})\ r}{(2N-1)u_r(r)^2},$$

\noindent for every $r\in (0,r_0)$, which is the desired conclusion in the interval $(0,r_0)$ for $M_1=\left((1-2^{1-2N})/((2N-1)K)\right)^{1/2}$. To finish the proof it remains to show that $u_r(r)<0$ for every $0<r\leq 1$. Indeed, if $u_r(r')\geq 0$ for some $0<r'\leq 1$ then, from the nonnegativeness and the monotonicity of $-r^{N-1}u_r(r)$ in $(0,1]$, it is deduced that $-r^{N-1}u_r(r)=0$ for every $0<r\leq r'$. Hence $u$ is constant in $(0,r']$, a contradiction. \qed

\section{Semi-stable radial weak solutions in a ball}

The following lemma gives a characterization of radial weak solutions of (\ref{theball}) and will be useful to prove Theorem \ref{N=2}.

\begin{lemma}\label{characterization}
Let $\Omega=B_1$, $f\in C(\mathbb{R})$ and $u$ be a radial function in $\overline{B_1}$. Then $u$ is a weak solution of (\ref{theball}) if and only if the following holds:
\begin{enumerate}
\item[(i)] $u\in C^2(0,1]$, $u(1)=0$ and $-\Delta u(x)=f(u(x))$ pointwise in $\overline{B_1}\setminus\{ 0\}$.
\item[(ii)] $f(u)\in L^1(B_1)$.
\item[(iii)] $\lim_{r\to 0} r^{N-1}u_r(r)=0$.
\end{enumerate}
\end{lemma}

\

\noindent\textbf{Proof.}
Let us prove first the necessary conditions. Suppose that $u$ is a radial weak solution of (\ref{theball}). Then it is well known that

$$u(r)=-\int_r^1 \left(\frac{u_r(1)+\int_t^1 s^{N-1}f(u(s))ds}{t^{N-1}}\right)dt,$$

\noindent and (i) is proved. On the other hand since $f(u)\delta\in L^1(B_1)$ then $f(u)\in L^1(B_{1/2})$. Taking into account that $f(u)$ is continuous in $\overline{B_1}\setminus B_{1/2}$, (ii) is proved. To prove (iii), consider $\zeta\in C^2(\overline{B_1})$ satisfying $\zeta=0$ on $\partial B_1$ and $\zeta=1$ in $B_{1/2}$. Applying (\ref{weakdefinition}) we deduce

$$0=\int_{B_1}\left(u\Delta \zeta+f(u)\zeta\right) dx=\lim_{r\to 0}\int_{B_1\setminus \overline{B_r}}\left(u\Delta \zeta+f(u)\zeta\right) dx$$ $$=\lim_{r\to 0}\int_{B_1\setminus \overline{B_r}}\left(u\Delta \zeta-\zeta\Delta u\right) dx=\lim_{r\to 0}\int_{\partial(B_1\setminus \overline{B_r})} \left(u\nabla \zeta-\zeta\nabla u\right)$$ $$=\lim_{r\to 0}\left(-\omega_N r^{N-1}u_r(r)\right),$$

\noindent and (iii) follows.

Suppose now that (i), (ii) and (iii) hold for a radial function $u$ defined in $\overline{B_1}$. From (iii) it is deduced that $\lim_{r\to 0}u(r)/\vert \log r\vert=0$ for $N=2$, while $\lim_{r\to 0}u(r)r^{N-2}=0$ for $N\geq 3$. In all the cases we have $\lim_{r\to 0}r^{N-1}u(r)=0$, which gives $r^{N-1}u(r)\in L^\infty(0,1)$ and then $u\in L^1(B_1)$. On the other hand (ii) clearly implies $f(u)\delta\in L^1(B_1)$. What is left to show is (\ref{weakdefinition}). To this end, consider $\zeta\in C^2(\overline{B_1})$ satisfying $\zeta=0$ on $\partial B_1$. Applying (i) and (ii) we obtain that

$$\int_{B_1}\left(u\Delta \zeta+f(u)\zeta\right) dx=\lim_{r\to 0}\int_{B_1\setminus \overline{B_r}}\left(u\Delta \zeta+f(u)\zeta\right) dx$$ $$=\lim_{r\to 0}\int_{B_1\setminus \overline{B_r}}\left(u\Delta \zeta-\zeta\Delta u\right) dx=\lim_{r\to 0}\int_{\partial(B_1\setminus \overline{B_r})} \left(u\nabla \zeta-\zeta\nabla u\right)=$$ $$\lim_{r\to 0}\int_{\partial B_r} \left(u\nabla \zeta-\zeta\nabla u\right).$$

Consider $M>0$ such that $\vert \zeta\vert$, $\vert \nabla \zeta \vert\leq M$ in $\overline{B_1}$. Applying $\lim_{r\to 0} r^{N-1}u(r)=0$ and
$\lim_{r\to 0} r^{N-1}u_r(r)=0$ the proof is complete by observing that

$$\left\vert\int_{\partial B_r} \left(u\nabla \zeta-\zeta\nabla u\right)\right\vert\leq \int_{\partial B_r} \left(\vert u\nabla \zeta\vert+\vert \zeta\nabla u\vert\right)\leq M\int_{\partial B_r} \left(\vert u\vert+\vert \nabla u\vert\right)$$ $$=M\omega_N r^{N-1}\left( \vert u(r)\vert+\vert u_r(r)\vert\right)\rightarrow 0 \mbox{ as } r\rightarrow 0. \qed $$

\begin{remark}
We can apply this characterization to the radial functions $u(r)=r^{-2/(p-1)}-1$, ($p>1$) mentioned in the Introduction. We have that $u$ is a solution of (\ref{mainequation}) for $f(u)=2\left(Np-2p-N\right)/(p-1)^2 (1+u)^p$. Applying Lemma \ref{characterization}, we check at once that $u$ is a radial weak solution of (\ref{theball}) if and only if $N\geq 3$ and $p>N/(N-2)$.
\end{remark}

\

\noindent\textbf{Proof of Theorem \ref{N=2}.}
Suppose that $u\notin H^1(B_1)$. Applying Theorem \ref{nonenergy} we have that there exist $M>0$ and $0<r_0<1$ such that $\vert u(r)\vert\geq M \vert\log r\vert$ for every $r\in (0,r_0)$. On the other hand, since $u$ is a radial weak solution of (\ref{theball}) we could apply (iii) of Lemma (\ref{characterization}) and obtain $\lim_{r\to 0}ru_r(r)=0$. In particular $\lim_{r\to 0}u(r)/\vert\log r\vert=0$, a contradiction.

Thus $u$ is an energy solution (i.e. $u\in H^1(B_1)$). It is known (see \cite{cabrecapella}) that $u\in L^\infty(B_1)$ and then, by standard regularity arguments, $u\in C^2(\overline{B_1})$. \qed

\end{document}